\magnification=1200
\documentstyle{amsppt}
\hsize=34 true pc
\vsize=53 true pc
\voffset=0.1 true in
\hoffset=1.0 true cm 
\def\fsl{\frak s\frak l({r+1},\text{\bf C})}  
\NoRunningHeads
\overfullrule0pt
 
\font\bigbold=cmbx10 scaled \magstep2 
\font\teni=cmmi10 scaled \magstep2
\font\seveni=cmmi7 scaled \magstep2
\font\fivei=cmmi5 scaled \magstep2
\textfont9=\teni \scriptfont9=\seveni \scriptscriptfont9=\fivei

\document
\baselineskip=16pt


\newfam\bigfam 
\topmatter
\title
\bigbold {The Polynomial Behavior of Weight Multiplicities}  
\bigbold {for the Affine Kac-Moody Algebras} $A^{(1)}_r$  
\endtitle
\medskip

\author
Georgia Benkart$^{*}$, 
Seok-Jin Kang$^{\dagger}$, 
Hyeonmi Lee$^{\dagger}$, \\
Kailash C. Misra$^{\ddagger}$,
Dong-Uy Shin$^{\dagger}$
\\
\\
\endauthor

\thanks 
\hskip -.2 truein {$^*$Supported in part by NSF Grants
$\#$DMS-9300523 and $\#$DMS-9622447}
\endthanks

\thanks 
\noindent {$^{\dagger}$Supported by the Non-directed Research Fund,
Korea  Research Foundation, 1996}
\endthanks

\thanks 
\noindent {$^{\ddagger}$Supported in part by NSA/MSP Grant
$\#$MDA904-96-1-0013
 }
\endthanks

\affil
\noindent 
$^{*}$Department of Mathematics \\
University of Wisconsin \\
Madison, Wisconsin 53706-1388, USA \\
benkart\@math.wisc.edu \\
\\
\noindent 
${\dagger}$Department of Mathematics \\
Seoul National University  \\
Seoul 151-742, Korea \\
sjkang\@math.snu.ac.kr \\
hmlee\@math.snu.ac.kr\\
dushin\@math.snu.ac.kr\\
\\
\noindent 
${\ddagger}$Department of Mathematics \\
North Carolina State University \\
Raleigh, North Carolina 27695-8205, USA \\
misra\@math.ncsu.edu \\
\endaffil

\abstract
We prove that the multiplicity of an arbitrary 
dominant weight for an integrable highest weight representation 
of the affine Kac-Moody algebra $A_{r}^{(1)}$ is a polynomial
in the rank $r$. In the process we show that the degree of this polynomial 
is less than or equal to the depth of the weight with respect to 
the highest weight. These results allow weight multiplicity
information for small ranks to be transferred to arbitrary ranks. 

\endabstract

\endtopmatter

\centerline {\bf Introduction}  

\vskip 4mm
\hskip 5mm The irreducible highest weight representations of affine Kac-Moody 
algebras have played an increasingly important
role in diverse areas of mathematics and physics. When its level is positive,
such a representation is infinite-dimensional.  It is 
parameterized by a dominant integral highest weight  and has finite-dimensional weight spaces.
The formal character of such a representation records the multiplicity of each
weight, and the well-known {\it Weyl-Kac character  formula} ([K2, p. 173])
provides a precise expression for the character.  The character formula involves
a sum over the Weyl group in both its numerator and denominator which makes it
impractical for explicitly computing multiplicities.  However, when the
character formula is applied to the one-dimensional trivial representation, it
gives the {\it denominator identity}, and from the denominator identity Peterson
[P] has derived Freudenthal-type  recursive formulas for calculating root and
weight multiplicities.   These formulas enabled Kass, Moody, Patera, and
Slansky [KMPS] to develop tables of weight multiplicities for certain weights of low
level irreducible highest weight representations for affine Kac-Moody algebras having rank
less than 8.  

\vskip 2mm
\hskip 5mm In 1987 while analyzing the weight multiplicities of the 
irreducible highest weight representations of the untwisted classical affine 
Kac-Moody algebras, Benkart and Kass (see [BK]) conjectured certain
polynomial behavior for the weight multiplicities of these representations
and introduced the notion of a ``rank-zero string function''. 
The conjectures were confirmed in [BKM2] for any irreducible highest weight representation
of the affine Kac-Moody algebras $A_r^{(1)}$ for weights having
depth $\leq 2$. In fact, in [BKM2] the multiplicities of such weights were 
given by explicit polynomials whose coefficients involve Kostka numbers. 
However, it seems to be very difficult to extend the methods of [BKM2],
which were based on the root multiplicity formula for Kac-Moody 
algebras obtained in [Ka2] and the representation theory of $\fsl$,  
to prove the conjecture for arbitrary depths. 

\vskip 2mm
\hskip 5mm  In this paper, we adopt a completely different approach to 
prove that the multiplicity of an arbitrary dominant weight for
an irreducible highest weight representation of the affine Kac-Moody algebra  
$A_r^{(1)}$ is a polynomial in the rank $r$.  Although the precise degree of these
polynomials is not determined in this work, an upper bound is obtained for the degree, and
this upper bound  coincides with the degree conjectured by Benkart and Kass 
(see [BKM2], Conjecture A).  
\vskip 2mm
\hskip 5mm  Briefly, our argument proceeds as follows:   Let $L(\lambda)$ denote the
irreducible highest weight $A_r^{(1)}$-module with highest weight $\lambda = \sum_{i=0}^r
a_i\Lambda_i -  m\delta$, where $\Lambda_0,\Lambda_1,\cdots,\Lambda_r$ are the fundamental
weights and $\delta$ is the null root.   We consider the minimal graded Lie algebra ${\Cal
L}$ with local part $L(\lambda)\oplus A_r^{(1)}\oplus L^*(\lambda)$, where 
$L^*(\lambda)$ is the finite dual space of $L(\lambda)$ (see Section 2). Then
${\Cal L}$ is isomorphic to the indefinite Kac-Moody algebra  
$\widehat{\frak g}$
associated to the Cartan matrix 
$\widehat {\frak A}=(a_{i,j})_{i,j=-1,0,1,\cdots,r}$,
whose first column consists of the entries
$2$, $-a_0$, $-a_1$, $\cdots$, $-a_r$. When 
the first row and the first 
column of $\widehat {\frak A}$ are deleted, the result is the Cartan matrix of the affine
Kac-Moody algebra $A_r^{(1)}$. Now
any weight $\mu$ of $L(\lambda)$ can be viewed as a root in 
$\widehat{\frak g}$, and its multiplicity as a root of  $\widehat{\frak g}$
is the same as its multiplicity as a weight of $L(\lambda)$.  
We use Peterson's recursive root multiplicity formula in conjunction with a tricky inductive
argument to establish the polynomial behavior of the dominant weights of
$L(\lambda)$.  
\vskip 2mm
\hskip 5mm The proof of the polynomial conjecture permits
much of the information in [KMPS] to be extended to arbitrary ranks, and 
it provides a means of relating string functions for various algebras. 
A better understanding of the polynomial nature of the multiplicities for more
general sequences of Kac-Moody algebras (beyond the $A_r^{(1)}$ case treated
here) would allow results about multiplicities for affine algebras to be
transferred to hyperbolic and indefinite Kac-Moody algebras, where only very
limited information is currently known.   

\vskip 5mm

\centerline {\bf Acknowledgments}

\hskip 5mm Part of this work was completed while Seok-Jin Kang visited
Yale University in February 1995 and University of Wisconsin-Madison
in the summer of 1995. He would like to express his sincere gratitude
to the Mathematics Departments of Yale University and the
University of Wisconsin-Madison for their hospitality.
We are very grateful to George Seligman 
and Hye-Kyung Oh for their interest in this work and many valuable
discussions. 
\vskip 8mm

\centerline {\bf \S 1. The affine weight lattice and the conjecture \rm}
\vskip 4mm
\hskip 5mm  Suppose that $I=\{ 0,1,\cdots, r \}$, and let ${\frak A}=(a_{i,j})_{i,j\in I}$
be the affine Cartan matrix of type $A_r^{(1)}$:
$${\frak A}= \pmatrix 
\hfill 2&\hfill -1&\hfill 0&\hfill 0&\hfill \cdots&\hfill 0 &\hfill -1 \,\, \\ 
\hfill-1&\hfill 2 &\hfill -1&\hfill 0&\hfill \cdots&\hfill 0 &\hfill 0 \,\, \\ 
\hfill 0&\hfill -1 &\hfill 2&\hfill -1&\hfill \cdots&\hfill 0 &\hfill 0 \,\, \\ 
\hfill 0&\hfill 0 &\hfill -1&\hfill 2&\hfill \cdots&\hfill 0 &\hfill 0 \,\, \\ 
\hfill \vdots&\hfill \vdots&\hfill \vdots&\hfill \vdots&\hfill \ddots&\hfill \vdots&\hfill \vdots \,\, \\  
\hfill 0&\hfill 0&\hfill 0&\hfill 0&\hfill \cdots&\hfill2 &\hfill -1 \,\, \\ 
\hfill -1&\hfill 0&\hfill 0&\hfill 0&\hfill \cdots&\hfill -1&\hfill 2 \,\, \\ 
\endpmatrix. \tag 1.1$$
Let ${\frak h}$ be a vector space over {\bf C} with
a basis $\{ h_0, h_1, \cdots, h_r, d \}$. Define linear functionals 
$\alpha_i \in {\frak h}^*$ ($i\in I$) by
$$ \alpha_i(h_j)=a_{j,i} \ \  \text {for} \ j\in I, \ \ 
\alpha_i(d)=\delta_{i,0}. \tag 1.2 $$
Then the triple $({\frak h}, \Pi=\{\alpha_i|\ i\in I \},
\Pi^{\vee}=\{h_i|\  i\in I \})$ provides a realization of the matrix ${\frak A}$ in the sense of
[K2, Chap.\,1].
The Kac-Moody algebra ${\frak g}$ associated with the affine matrix ${\frak A}$
is the {\it affine Kac-Moody algebra of type $A_r^{(1)}$}. 
We denote by $e_i$, $f_i$, $h_i$ $(i\in I)$ and $d$ the generators of the 
algebra ${\frak g}$. The subalgebra ${\frak g}_0$ of ${\frak g}$ generated 
by $e_i$, $f_i$, $h_i$ ($i=1, \cdots, r$) is a finite-dimensional simple
Lie algebra of type $A_r$ which is 
isomorphic to the Lie algebra $\fsl$ of $(r+1)\times
(r+1)$ complex matrices of trace zero. 

\vskip 2mm
\hskip 5mm Let $c=h_0+h_1+ \cdots +h_r$. Then $[c,x]=0$ for all $x\in 
{\frak g}$, and $c$ is the {\it canonical central element} of
${\frak g}$. Note that $\{ h_1, \cdots, h_r, c, d \}$ forms another 
basis of ${\frak h}$. Since the matrix ${\frak A}$ is symmetric, there is a 
nondegenerate symmetric bilinear form on ${\frak h}$ which satisfies
$$
\aligned
& (h_i|h_j)=a_{i,j} \ \ \text {for} \ i,j=1, \cdots, r, \\
& (h_i|c)=(h_i|d)=0 \ \ \text {for} \ i=1,\cdots, r, \\
& (c|c)=(d|d)=0, \ \ (c|d)=1.
\endaligned
\tag 1.3
$$
Define linear functionals $\Lambda_i \in {\frak h}^*$ $(i\in I)$ and 
$\delta \in {\frak h}^*$  by
$$
\aligned
& \Lambda_i (h_j) =\delta_{i,j}, \ \ \Lambda_i(d)=0, \\
& \delta(h_j)=0, \ \ \delta (d)=1 \ \ \text {for} \ j\in I.
\endaligned \tag 1.4
$$
Then $\delta$ can be expressed as 
$\delta=\alpha_0+\alpha_1 + \cdots + \alpha_r.$
It is easy to see that $\{ \Lambda_0, \Lambda_1, \cdots, \Lambda_r, \delta \}$
and $\{\Lambda_0, \alpha_0, \alpha_1, \cdots, \alpha_r\}$ are both bases of
the complex vector space ${\frak h}^*$, and
$$
\aligned
& \alpha_i = -\Lambda_{i-1} + 2\Lambda_i-\Lambda_{i+1} + \delta_{i,0}\delta \ \ 
(i \; \text {mod} \, r+1). \endaligned \tag 1.5
$$   

\vskip 2mm
\hskip 5mm The {\it affine weight lattice} $P$ is defined to be
$P=\text {\bf Z}\Lambda_0 \oplus \text {\bf Z}\Lambda_1 \oplus 
\cdots \oplus \text {\bf Z}\Lambda_r \oplus \text {\bf Z}\delta,$
and the elements of
$P^+=\{\lambda\in P| \ \lambda(h_i) \in \text {\bf Z}_{\ge 0} \  \text 
{for all} \  i\in I \}$ are the {\it dominant integral weights} for 
the algebra ${\frak g}$. 
For $\lambda, \mu \in P$, we say that $\mu$ is {\it related} to $\lambda$,
which we denote by $\mu \sim \lambda$, if $\lambda-\mu \in Q$, where 
$Q = \bigoplus_{i = 0}^r \text {\bf Z}\alpha_i$ is the root lattice. 
For example,
if $\mu \in P$ is a weight of the irreducible highest weight module
$L(\lambda)$ over ${\frak g}$ with highest weight $\lambda\in P$,
then $\lambda-\mu \in Q_{+} = \bigoplus_{i = 0}^r \text {\bf Z}_{\ge 0}\alpha_i$, and
hence $\mu$ is related to $\lambda$.

\vskip 2mm
\hskip 5mm Let $l$ be a positive integer. 
A dominant integral weight $\lambda \in P^+$  is said to have 
{\it level} $l>0$ if
$\lambda(c)=l$. The weight $\lambda$ can be uniquely expressed in the form 
$$\lambda=a_0 \Lambda_0 + a_1 \Lambda_1 + \cdots + a_r \Lambda_r 
-m\delta, \tag 1.6$$
where $m\in \text {\bf Z}$ and $a_i \in \text {\bf Z}_{\ge 0}$ for
$i=0,1,\cdots, r$.
Since $c=h_0+h_1+\cdots +h_r$, we have
$$\lambda(c)=a_0+a_1+\cdots +a_r=l. \tag 1.7$$

\hskip 5mm Let $\mu \in P$ be an integral weight and suppose $\mu$ is
related to $\lambda$. We write 
$$\mu=b_0 \Lambda_0 + b_1 \Lambda_1 + \cdots + b_r \Lambda_r 
-n\delta, \tag 1.8$$
where $n \in \text {\bf Z}$ and $b_i \in  \text {\bf Z}$ for $i=0,1,
\cdots, r$, and we set $d_i=b_i-a_i$ $(i\in I)$. 
Since $\mu$ is related to $\lambda$, we can write
$\mu=\lambda-\sum_{i=0}^r k_i \alpha_i$ for some $k_i \in \text {\bf Z}$
$(i=0, 1, \cdots, r)$. Therefore, since $\alpha_i(c) = 0$ for all $i$,  we must
have $$\mu(c)=b_0+b_1+\cdots +b_r=l, $$
which implies 
$$d_0+d_1+\cdots +d_r=0. \tag 1.9$$

\hskip 5mm Using the linear system
$$
\aligned
& \mu(h_j)=a_j-\sum_{i=0}^r k_i a_{j,i} =b_j  
\ \ \text {for} \ j=0,1,\cdots,r,\\
& \mu(d)=-m-k_0=-n,
\endaligned
\tag 1.10
$$
we can solve for the $k_i$'s to obtain 
$$
\aligned
 k_0&=n-m, \\
 k_i&=n-m+d_{i+1}+2d_{i+2}+\cdots +(r-i)d_r-(r-i+1)\dsize\frac {N}{r+1} \\
    &\hskip .9 truein  \text {for} \
i=1,\cdots, r, \endaligned
\tag 1.11
$$
where $N=d_1+2d_2+\cdots +rd_r$. In particular,
$$k_r=n-m-\dsize\frac {N}{r+1} \in \text {\bf Z}.$$
Thus 
$$N =d_1+2d_2+ \cdots + rd_r \equiv 0 \ \ \text {mod} \  r+1. \tag 1.12$$

\hskip 5mm Conversely, suppose $\mu=\sum_{i=0}^r b_i \Lambda_i -n\delta
\in P$ is an integral weight with $b_i \in \text {\bf Z}$ $(i\in I)$
and $n\in \text {\bf Z}$ which satisfies (1.9) and (1.12) for
$d_i = b_i - a_i\ \ (i \in I)$.  
Then we can write $\mu=\lambda-\sum_{i=0}^r k_i \alpha_i$, where the
$k_i$ are given by (1.11). Hence $\mu$ is related to $\lambda$. 
Therefore, we obtain:

\vskip 4mm 
\proclaim{ Proposition 1.13} \ Let $\lambda=\sum_{i=0}^r a_i \Lambda_i
-m\delta$ $(a_i \in \text 
{\bf Z}_{\ge 0}, m\in \text {\bf Z})$ be a dominant integral weight of level $l$, and let  
$\mu=\sum_{i=0}^r b_i \Lambda_i -n\delta \in P$ $(b_i \in \text {\bf Z}, 
n \in \text {\bf Z})$ be an integral weight. 
Then $\mu$ is related to $\lambda$ if and only if 
$$
\aligned
& d_0+d_1+\cdots +d_r =0, \\
& d_1 +2d_2 + \cdots + rd_r \equiv 0 \ \ \text {mod} \ r+1,  
\endaligned
\tag 1.14
$$
\noindent where $d_i = b_i-a_i \ \ (i \in I)$.  
\endproclaim

\vskip 3mm 
\hskip 5mm  
Now fix a positive integer $l$ and a dominant integral weight $\lambda$
of level $l$ given by (1.6).  Since $\lambda(c)=a_0+a_1+\cdots+a_r=l$,
there must be a {\it gap} in the expression (1.6) for 
$\lambda$ if $r \ge l$. That is, there exist nonnegative integers $s$ and $t$
with $s+t \le r$ such that
$$a_{s-1} \neq 0, \ \ a_{r-t+1} \neq 0, \ \ 
a_{s}=a_{s+1}=\cdots =a_{r-t}=0. \tag 1.15$$

\hskip 5mm 
The viewpoint we adopt here is that the weight $\lambda$ is completely determined by
the following data: (i) \ an $s$-tuple of nonnegative integers 
$\underline a=(a_0, a_1, \cdots,$ 
$a_{s-1})$ with $a_{s-1} \neq 0$, \ (ii) 
a $t$-tuple of nonnegative integers 
$\underline a'=(a_{r-t+1}, a_{r-t+2}, \cdots,$
$a_{r})$ with $a_{r-t+1} \neq 0$, and (iii) an integer $m$. 
Note that this determining data is independent of $r$.
Thus, dominant integral weights will be regarded as the same for all $r \ge l$  
provided they have the same determining data. It is important to observe that
a different choice of gap in expression (1.6) yields a 
different weight. For example, consider 
$\lambda=2\Lambda_2+\Lambda_4-\delta$ when $r=5$.
If we take $s=5$ and $t=0$, then 
the determining data for $\lambda$ is
$\underline a=(0,0,2,0,1)$, $\underline a'=\emptyset$, $m=1$,
and $\lambda$ can be written as 
$\lambda=2\Lambda_2+\Lambda_4-\delta$ for all $r\ge 5$.
On the other hand, we can choose a different gap by taking
$s=3$, $t=2$. In this case, the determining data for $\lambda$
is $\underline a=(0,0,2)$, $\underline a'=(1,0)$, $m=1$,
and $\lambda$ can be expressed as $\lambda=2\Lambda_2+\Lambda_{r-1}-\delta$
for all $r\ge 5$.   

\vskip 2mm
\hskip 5mm Suppose the dominant integral weight $\lambda$ is given by the
determining data 
$\underline a=(a_0, a_1, \cdots,$ 
$a_{s-1})$, 
$\underline a'=(a_{r-t+1}, a_{r-t+2}, \cdots,$
$a_{r})$, and $m$. 
Let $\mu=\sum_{i=0}^r b_i \Lambda_i -n \delta$ $(b_i \in \text {\bf Z}_{\ge 0},
n\in \text {\bf Z})$ be a dominant integral weight
of level $l$. Then  
$\mu(c)=b_0+b_1+\cdots+b_r=l$, and if $r\ge l$, there must be 
a gap in the expression (1.8) for $\mu$.
Moreover, if $r \ge l+s+t$, then the gap of $\lambda$ is sufficiently
large that there exists a gap of $\mu$ which overlaps the
gap of $\lambda$.  As a result, we can associate determining data
$\underline b=(b_0, b_1, \cdots, b_{s'-1})$,
$\underline b'=(b_{r-t'+1}, b_{r-t'+2}, \cdots, b_r)$,
and $n\in \text {\bf Z}$ to  $\mu$,  
where  $b_{s'-1} \neq 0$,
$b_{r-t'+1} \neq 0$, and  $s'$, $t'$ are nonnegative 
integers satisfying $s'+t'\le r$, $s+t'\le r$, and $s'+t\le r$. 
Hence, if we let $p=\text {max} (s,s')$ and 
$q=\text {max} (t,t')$, 
the weights $\lambda$ and $\mu$ share a 
{\it common gap}:
$$
\aligned
&a_{p}=a_{p+1}= \cdots = a_{r-q}=0, \\
&b_{p}=b_{p+1}=\cdots = b_{r-q}=0. 
\endaligned
\tag 1.16
$$
(Note that $p \le r-t$ and $s \le r-q$.)

\vskip 2mm
\hskip 5mm  From now on we assume that $r \ge l+s+t$. 
If $\mu$ is related to $\lambda$
for infinitely many values of $r\ge l+s+t$,  
then the congruence equation (1.14) holds 
for all those values of $r$.  Hence,  

$$
\aligned
N &= d_1 + 2d_2 + \cdots + rd_r \\
& =  \big(d_1+2d_2+\cdots
+(p-1)d_{p-1}\big) \\
&\quad \; -\big(qd_{r-q+1}+(q-1)d_{r-q+2}+\cdots+ 2d_{r-1}+d_r\big) \\
& \quad \quad \; + (r+1)\big(d_{r-q+1}+d_{r-q+2}+\cdots + d_r\big) \\
\endaligned
\tag 1.17
$$
is divisible by $r+1$ for all such values of $r$,  which implies  

$$ d_1+2d_2+\cdots + (p-1)d_{p-1} =qd_{r-q+1}+\cdots+ 2d_{r-1}+d_r.
\tag 1.18
$$

\hskip 5mm We now define 
$$
\aligned
d_{\lambda}(\mu) &= n-m-(d_1+2d_2+\cdots + (p-1)d_{p-1}) \\
&= n-m -(qd_{r-q+1}+\cdots +2d_{r-1}+d_r),
\endaligned
\tag 1.19
$$
and refer to $d_{\lambda}(\mu)$ as the {\it depth of $\mu$ 
with respect to $\lambda$}.
It follows from (1.11) and (1.19) that 
$$
k_i = \left \{\aligned & d_{i+1}+2d_{i+2}+\cdots + (p-i-1) d_{p-1} 
+d_{\lambda}(\mu) \\
& \ \ \ \ \ \ \ \ \ \  \text {for} \ i = 0,1,\cdots, p-2, \\
& d_{\lambda} (\mu) \ \ \text {for} \ i=p-1, p,\cdots, r-q, r-q+1, \\
& (i-(r-q+1))d_{r-q+1}+\cdots +2d_{i-2}+ d_{i-1} +d_{\lambda}(\mu)\\
& \ \ \ \ \ \ \ \ \ \ \text {for} \ i= r-q +2,\cdots, r-1, r.
\endaligned \right. 
\tag 1.20
$$ 

\vskip 2mm
\hskip 5mm For $i=0,1,\cdots, r$, let $m_i=k_i-d_{\lambda}(\mu)$,
and define $\mu_0=\lambda-\sum_{i=0}^r m_i \alpha_i$. Then 
$\mu=\mu_0-d_{\lambda}(\mu) \delta$ and
$d_{\lambda} (\mu_0)=0$.
Note that $\mu(h_i)=\mu_0(h_i)=b_i$ for all $i=0,1,\cdots,r$.
By (1.20), we obtain
$$
m_i = \left \{\aligned
& d_{i+1}+2d_{i+2}+\cdots + (p-i-1) d_{p-1} \\
& \ \ \ \ \ \ \ \ \ \ \text {for} \ i = 0,1,\cdots, p-2, \\
& 0 \ \ \text {for} \ i=p-1, p, \cdots, r-q, r-q+1, \\
& (i-(r-q+1))d_{r-q+1}+\cdots +2d_{i-2}+ d_{i-1} \\
& \ \ \ \ \ \ \ \ \ \  \text {for} \ i= r-q +2,\cdots, r-1, r.
\endaligned \right.
\tag 1.21
$$
In particular, the values of $m_i$'s do not depend on $r$. 

\vskip 2mm 
To summarize the above discussion, we have

\vskip 2mm

\proclaim{ Proposition 1.22} 
\ Let $\lambda \in P^+$ 
be a dominant integral weight of level $l>0$ with
determining data 
$\underline a=(a_0, a_1, \cdots, a_{s-1})$,
$\underline a'=(a_{r-t+1}, a_{r-t+2}, \cdots, a_{r})$, 
and $m \in \text {\bf Z}$.
Assume that $r\ge l+s+t$ and that 
$\mu$ is a dominant integral weight
of level $l$ with
determining data 
$\underline b=(b_0, b_1, \cdots, b_{s'-1})$,
$\underline b'=(b_{r-t'+1}, b_{r-t'+2}, \cdots, b_{r})$, 
and $n \in \text {\bf Z}$ such that $s'+t' \le r$, $s+t' \le r$, and $s'+t \le r$.
If $\mu$ is related to $\lambda$ for infinitely 
many values of  $r\ge l+s+t$,
then $\mu$ can be uniquely written as 
$\mu=\mu_0-d_{\lambda}(\mu)\delta$, 
where $d_{\lambda}(\mu)$ and
$\mu_0=\lambda-\sum_{i=0}^r m_i\alpha_i$  are determined by (1.19) and (1.21) 
for $d_i = b_i-a_i \ 
(i \in I)$, and  $d_{\lambda}(\mu_0)=0$. 
\endproclaim
\vskip 2mm
The following  lemma plays an important role in proving our main theorem (Theorem 3.4).

\vskip 2mm

\proclaim{ Lemma 1.23} 
\ Let $\lambda\in P^+$
be a dominant integral weight of level $l>0$ 
with determining data
$\underline a=(a_0, a_1, \cdots, a_{s-1})$,
$\underline a'=(a_{r-t+1}, a_{r-t+2}, \cdots, a_{r})$,
and $m \in \text {\bf Z}$.
Assume that $r \ge l+s+t$. Let $\mu \in P^+$ be a dominant 
integral weight of level $l$ with determining data
$\underline b=(b_0, b_1, \cdots, b_{s'-1})$,
$\underline b'=( b_{r-t'+1}, b_{r-t'+2}, \cdots, b_{r})$, and $n \in \text {\bf Z}$
such that $s'+t' \le r,\;  s+t' \le r$, and $s'+t \le r$.
Let $\tau \in P^+$ be a dominant 
integral weight of level $l$ with determining data
$\underline c=(c_0, c_1, \cdots, c_{s''-1})$,
$\underline c'=(c_{r-t''+1}, c_{r-t''+2}, \cdots, c_{r})$,
$n' \in \text {\bf Z}$
satisfying 
$s''+t'' \le r,\;  s+t'' \le r$, and $s''+t \le r$.
If $\mu \le \tau \le \lambda$, then $d_{\lambda}(\tau) \le d_{\lambda}(\mu)$.
\endproclaim

\noindent    
{\it Proof.} \
By Proposition 1.22, $\mu$ and $\tau$ can be expressed as follows:
$$\mu=\lambda-\sum_{i=0} ^r m_i\alpha_i-d_{\lambda}(\mu)\delta,$$ 
$$\tau=\lambda-\sum_{i=0} ^r m'_i\alpha_i-d_{\lambda}(\tau)\delta.$$ 
The condition $\tau \ge \mu$ implies

$$\tau-\mu=\sum_{i=0}^r (m_i-m'_i) \alpha_i
+(d_{\lambda}(\mu)-d_{\lambda}(\tau))\delta \in Q_{+}.
\tag 1.24$$
Since $a_j=0$ for all $j=s,s+1, \cdots, r-t$, we have
$$d'_j=c_j-a_j=c_j \ge 0 \ \ \text {for all} \ j=s, s+1, \cdots, 
r-t. \tag 1.25$$

Let $p= \text{max}(s,s'), \,q=\text{max}(t,t'), \,x= \text{max}(s,s'')$, and 
$y=\text{max}(t,t'')$.  Recall that $x \le r-t$ and $s \le r-y$, and assume that
$i\in \{s-1, s,\cdots, r-t, r-t+1 \}$.   
It suffices to consider the following three cases: 
$$ \text {(i)} \ s-1 \le i \le x-2, \ \ \text {(ii)} \ x-1 \le i \le r-y+1, 
\ \ \text {(iii)} \ r-y+2 \le i \le r-t+1.$$
If $s-1\le i \le x-2$, then since $x \le r-t$, it follows from 
(1.21) and (1.25) that 
$$m_i'=d'_{i+1}+2d'_{i+2}+ \cdots + (x-i-1) d'_{x-1} \ge 0.$$
If $x-1 \le i \le r-y+1$, then by (1.21) $m'_i=0$.
Finally, if $r-y+2 \le i \le r-t+1$, then because $s \le r-y$,
(1.21) and (1.25) yield 
$$m'_i=(i-(r-y+1)) d'_{r-y+1}+ \cdots +2 d'_{i-2}+d'_{i-1} \ge 0.$$
Therefore, $m'_i \ge 0$ for all $i=s-1,s, \cdots, r-t, r-t+1$.
In particular, since $p \ge s$, $q \ge t$, we conclude that
$m'_i \ge 0$ for all $i=p-1,p, \cdots, r-q, r-q+1$.

\vskip 2mm
\hskip 5mm Observe that in (1.24)
the coefficient of $\alpha_i$ in $\tau-\mu$ for $i=p-1,p, \cdots, r-q,
r-q+1$ is
$$-m'_i+ d_{\lambda}(\mu)-d_{\lambda}(\tau)  \ge 0.  \tag 1.26$$ 
Since $m'_i \ge 0$ for $i=p-1,p, \cdots, r-q, r-q+1$, 
it must be that $d_{\lambda}(\tau) \le d_{\lambda}(\mu)$.\hskip 1cm $\square$

\vskip 2mm
\hskip 5mm The following is a more detailed formulation of a conjecture
presented in [BK].

\vskip 2mm
\noindent
\proclaim {\bf Conjecture} \ Let $\lambda \in P^+$ 
be a dominant integral weight of level $l>0$ with
determining data 
$\underline a=(a_0, a_1, \cdots, a_{s-1})$,
$\underline a'=(a_{r-t+1}, a_{r-t+2}, \cdots, a_{r})$, and
$m \in \text {\bf Z}$.
Assume that $r\ge l+s+t+2$, and let 
$L(\lambda)$ be the irreducible highest weight 
module over the affine Kac-Moody algebra ${\frak g}$ of type $A_r^{(1)}$ with
highest weight $\lambda$. 
Let $\mu \in P^+$ be a dominant
integral weight of level $l$ with determining data 
$\underline b=(b_0, b_1, \cdots, b_{s'-1})$,
$\underline b'=(b_{r-t'+1}, b_{r-t'+2}, \cdots, b_{r})$, and
$n \in \text {\bf Z}$ such that $s'+t' \le r$,\, $s+t' \le r$, and $s'+t \le r$.
Suppose that $\mu$ is related to $\lambda$ for 
infinitely many values of $r\ge l+s+t+2$.
If $\mu$ is a weight of the ${\frak g}$-module $L(\lambda)$ for some
$r_0 \ge l+s+t+2$, then it is a weight of $L(\lambda)$ for all 
$r \ge r_0$, and the multiplicity of  $\mu$ 
in $L(\lambda)$ is given by a polynomial in $r$ of degree $d_{\lambda}(\mu)$.
If $d_{\lambda}(\mu) <0$, then the multiplicity of $\mu$ in $L(\lambda)$
is zero. \endproclaim

\vskip 2mm
\hskip 5mm In Section 3, we will prove a weaker version of the above 
conjecture. We will show that the multiplicity of $\mu$ in $L(\lambda)$
is given by a polynomial in $r$ of degree $\le d_{\lambda}(\mu)$. 
Our approach is to apply Peterson's formula to a certain indefinite 
Kac-Moody algebra ${\Cal L}$ which will be constructed in the next section.

\vskip 8mm

\centerline {\bf \S 2. The indefinite Kac-Moody algebra ${\Cal L}$ \rm}
\vskip 4mm
\hskip 5mm Recall that the Cartan subalgebra ${\frak h}$ of the affine 
Kac-Moody algebra ${\frak g}$ of type $A_r^{(1)}$ has the basis $\{h_1,\cdots,
h_r, c, d\}$, where $c=h_0+h_1+\cdots +h_r$ is the canonical central element, and
there is a nondegenerate symmetric bilinear form $(\ |\ )$ on ${\frak g}$
whose values on ${\frak h}$ are given by (1.3). Since the form $(\ |\ )$ is
nondegenerate on ${\frak h}$, for every $\mu \in {\frak h}^*$ there is a
unique element $t_{\mu}$ in ${\frak h}$ such that $\mu(h)=(h| t_{\mu})$ for
all $h\in  {\frak h}$. Thus the form $(\ |\ )$ induces a nondegenerate
symmetric  bilinear form on ${\frak h}^*$, also denoted by $(\ |\ )$, defined
by $(\mu | \nu )=(t_{\mu} | t_{\nu})$ for all $\mu, \nu \in {\frak h}^*$. 
In particular, $t_\delta = c$.   We
take a basis for $\frak h$ and extend it to a basis $\{x_i|\ i\in \Omega \}$
of $\frak g$ by adding     basis elements for each of the root spaces ${\frak
g}_\alpha$.   Since $({\frak g}_\alpha | {\frak g}_\beta) = 0$  unless $\beta
= -\alpha$, the dual basis $\{y_i|\ i\in \Omega \}$ of $\frak g$ with respect
to the form $(\ |\ )$ also consists of vectors in $\frak h$ and root
vectors.   

\vskip 2mm
\hskip 5mm Assume $\lambda=\sum_{i=0}^r a_i \Lambda_i -m\delta$ is a dominant 
integral weight of level $l>0$ for ${\frak g}$
as in Section 1. 
Let $L(\lambda)$ be the irreducible highest weight 
${\frak g}$-module with highest weight $\lambda$. The finite dual space  
$L^*(\lambda)$ is the irreducible lowest weight ${\frak g}$-module 
with lowest weight $-\lambda$, where the 
${\frak g}$-module action is given by 
$$\langle g \cdot v^*, w \rangle = -\langle v^*, g\cdot w \rangle 
\tag 2.1$$
for $g\in {\frak g}$, $v^*\in L^*(\lambda)$, $w\in L(\lambda)$ (see [K2, p. 149]).
Define a linear map $\psi : L^*(\lambda) \otimes L(\lambda)
\rightarrow  {\frak g}$ by
$$ \psi (v^*\otimes w)=- \dsize\frac {2}{(\lambda|\lambda)}
\sum_{i\in \Omega} \langle v^*, x_i \cdot w\rangle y_i,
\tag 2.2$$
where $\{x_i|\ i\in \Omega \}$ and $\{y_i|\ i\in \Omega \}$ are dual
bases of ${\frak g}$ as above. 
Then $\psi$ is a well-defined ${\frak g}$-module homomorphism, (compare with
[FF], [Ka1], [BKM1]), and hence 
the space $L(\lambda)\oplus {\frak g} \oplus L^*(\lambda)$
has the structure of a local Lie algebra with the bracket defined by
$$
\aligned
& [v^*, w]=\psi (v^*\otimes w), \\
& [g,w]=g\cdot w, \ \ [g,v^*]=g\cdot v^*
\endaligned
\tag 2.3
$$
for $g\in {\frak g}$, $v^*\in L^*(\lambda)$, $w \in L(\lambda)$ (see [K1]).

\vskip 2mm
\hskip 5mm Let ${\Cal F}_{+}$ (resp. ${\Cal F}_{-}$) be the free Lie
algebra generated by $L^*(\lambda)$ (resp. $L(\lambda)$),
and for $k\ge 1$, let ${\Cal F}_k$ (resp. ${\Cal F}_{-k}$) be the 
subspace of ${\Cal F}_{+}$ (resp. ${\Cal F}_{-}$) spanned by the vectors
of the form $[u_1 [u_2 [ \cdots [u_{k-1}, u_k] \cdots ]]]$ with
$u_j \in L^*(\lambda)$ (resp. $L(\lambda)$). 
In particular, ${\Cal F}_1=L^*(\lambda)$ and 
${\Cal F}_{-1}=L(\lambda)$. Let ${\Cal F}_0={\frak g}$ and
define 
$${\Cal F}={\Cal F}_{-}\oplus {\Cal F}_0 \oplus {\Cal F}_{+}
= \bigoplus_{k\in \text {\bf Z}} {\Cal F}_{k}.$$
Then ${\Cal F}$ is the maximal graded Lie algebra with local part
$L(\lambda) \oplus {\frak g} \oplus L^*(\lambda)$. 

\vskip 2mm
\hskip 5mm For $k\ge 2$, define the subspaces ${\Cal J}_{\pm k}$ 
of ${\Cal F}_{\pm k}$ by
$${\Cal J}_{\pm k}=\{ v\in {\Cal F}_{\pm k}| \ 
[u_1 [u_2 [ \cdots [ u_{k-1}, v] \cdots ]]]=0 \ \ 
\text {for all} \ u_i \in {\Cal F}_{\mp 1} \}, \tag 2.4$$
and let ${\Cal J}_{\pm} =\bigoplus_{k\ge 2} {\Cal J}_{\pm k}$.
Then ${\Cal J}_{\pm}$ is a graded ideal of ${\Cal F}_{\pm}$, and
${\Cal J}={\Cal J}_{-} \oplus {\Cal J}_{+}$ is the maximal graded ideal 
of ${\Cal F}$ which intersects the local part 
$L(\lambda) \oplus {\frak g} \oplus L^*(\lambda)$
trivially ([K1], [FF], [Ka1], [BKM1]). The Lie algebra 
${\Cal L}={\Cal F}/{\Cal J}=\bigoplus_{k\in \text {\bf Z}} {\Cal L}_{k}$
is the minimal graded Lie algebra with local part 
$L(\lambda) \oplus {\frak g} \oplus L^*(\lambda)$,
where ${\Cal L}_k={\Cal F}_k / {\Cal J}_k$ for $k\in \text {\bf Z}$.
In particular, ${\Cal L}_{-1}={\Cal F}_{-1}=L(\lambda)$,
${\Cal L}_0={\frak g}$, and ${\Cal L}_1={\Cal F}_1=L^*(\lambda)$.

\vskip 2mm
\hskip 5mm  Alternately,  let $\alpha_{-1}=-\lambda$ and consider
the Cartan matrix $\widehat{\frak A} =(a_{i,j})$ $(i,j=-1,0,1,\cdots,r)$
given by
$$a_{i,j}=\dsize\frac {2(\alpha_i|\alpha_j)}{(\alpha_i|\alpha_i)} \ \ 
\text {for} \ i,j=-1,0,1,\cdots,r. \tag 2.5$$
The first column of the matrix $\widehat {\frak A}$ consists of the entries
$2$, $-a_0$, $-a_1$, $\cdots$, $-a_r$, 
and deleting the first row and the first 
column of $\widehat {\frak A}$ gives the affine Cartan matrix of type 
$A_r^{(1)}$. If we let 
$$h_{-1}=-\dsize\frac {2 t_{\lambda}}{(\lambda|\lambda)},$$
where $t_{\lambda}\in {\frak h}$ is such that
$\lambda(h)=(h|t_{\lambda})$ for all $h\in {\frak h}$, then the triple
$( {\frak h}, \Pi=\{\alpha_{-1}, \alpha_0, \alpha_1, \cdots, \alpha_r \},
\Pi^{\vee}=\{ h_{-1}, h_0, h_1, \cdots, h_r \})$ provides a realization of the 
matrix $\widehat {\frak A}$. Let $\widehat {\frak g}$ be the indefinite Kac-Moody 
algebra associated with the Cartan matrix $\widehat {\frak A}$. 
It is a direct consequence of the Gabber-Kac theorem (see [GK]) that the following holds:

\vskip 2mm
\proclaim{ Proposition 2.6} {\rm ([FF], [Ka1], [BKM1])} \ 
Let $v_0$ (resp. $v_0^*$)
be the highest (resp. lowest)  weight vector of $L(\lambda)$
(resp. $L^*(\lambda)$) such that $\langle v_0^*, v_0\rangle =1$.
If $\widehat {\frak g}$ is the indefinite Kac-Moody algebra with the Cartan 
matrix $\widehat {\frak A}$ given by (2.5), then there is an isomorphism of
Lie algebras $\widehat {\frak g} \cong {\Cal L}$ defined by
$$
\aligned
& e_i \mapsto e_i, \ \ f_i \mapsto f_i, \ \ h_i \mapsto h_i 
\ \ \text {for} \ i=0,1,\cdots,r,\\
& e_{-1} \mapsto v_0^*, \ \ f_{-1} \mapsto v_0, \ \ 
 h_{-1} \mapsto - \dsize\frac {2 t_{\lambda}}{(\lambda|\lambda)}.
\endaligned
\tag 2.7
$$
\endproclaim

\vskip 2mm
\hskip 5mm It follows from Proposition 2.6 that the subspace
 ${\Cal L}_{\pm k}$ is the sum of all the root spaces $\widehat {\frak g}_{\pm \alpha}$, where
$\alpha$ is of the form $\pm (k\alpha_{-1}+\sum_{i=0}^r k_i \alpha_i)$ with
$k, k_i \in \text {\bf Z}_{\ge 0}$. 
In particular, the roots of the form
$\pm (\sum_{i=0}^r k_i \alpha_i)$ are roots of the affine Kac-Moody algebra
${\frak g}$, and the roots of the form $-\alpha_{-1}-\sum_{i=0}^r 
k_i \alpha_i$ are weights of the irreducible highest weight ${\frak g}$-module
$L(\lambda)$. Thus to compute the weight multiplicity of
$\lambda-\sum_{i=0}^r k_i \alpha_i$ in $L(\lambda)$, it suffices
to compute the root multiplicity of $-\alpha_{-1}-\sum_{i=0}^r k_i \alpha_i$
in the indefinite Kac-Moody algebra ${\Cal L}\cong \widehat {\frak g}$.

\vskip 8mm

\centerline {\bf \S 3. The weight multiplicity polynomials \rm}
\vskip 4mm
\hskip 5mm In this section, we will prove our main result. 
Fix a positive integer $l$ and 
a dominant integral weight $\lambda$ of level $l$ 
with determining data 
$\underline a=(a_0, a_1, \cdots, a_{s-1})$,
$\underline a'=(a_{r-t+1}, a_{r-t+2}, \cdots, a_{r})$, and
$m \in \text {\bf Z}$.
Assume that $r\ge l+s+t+2$, and 
let $L(\lambda)$ be the irreducible highest weight module over
the affine Kac-Moody algebra ${\frak g}=A_r^{(1)}$ with highest
weight $\lambda$.
Suppose $\mu$ is a dominant integral weight  of level $l$ with determining
data 
$\underline b=(b_0, b_1, \cdots, b_{s'-1})$,
$\underline b'=(b_{r-t'+1}, b_{r-t'+2}, \cdots, b_{r})$,
and $n \in \text {\bf Z}$, where $s'$, $t'$ are nonnegative integers 
satisfying $s'+t'\le r$, \, $s+t'\le r$, and $s'+t\le r$. 
Since the determining data associated to $\lambda$ and $\mu$ is fixed, the
integers $s, t, s', t'$ are all fixed. In particular, the integers 
$p=\text {max} (s,s')$ and $q=\text {max} (t,t')$ are fixed also.
Suppose further that $\mu$ is related to $\lambda$ for
infinitely many values of $r\ge l+s+t+2$. 
Our aim is to prove that if $\mu$ is a weight of $L(\lambda)$
for some $r_0 \ge l+s+t+2$, then it is a weight of $L(\lambda)$ for
all $r\ge r_0$, and the multiplicity of $\mu$ in $L(\lambda)$
is a polynomial in $r$ of degree $\leq d_\lambda(\mu)$, the depth of 
$\mu$ with respect to $\lambda$.  

\vskip 2mm
\hskip 5mm 
Let ${\Cal L} \cong \widehat {\frak g}$ be the minimal graded Lie algebra
with local part $L(\lambda)\oplus {\frak g} \oplus L^*(\lambda)$ 
constructed in Section 2.
Let $\widehat Q = \bigoplus_{i = -1}^r \text {\bf Z}\alpha_i$ denote 
the root lattice of
$\widehat {\frak g}$ with respect to the Cartan subalgebra $\frak h$.    
The roots of $\widehat {\frak g}$ belong to $\widehat Q_+ \cup \widehat Q_-$, where
$\widehat Q_+  = \bigoplus_{i = -1}^r \text {\bf Z}_{\geq 0}\alpha_i = -\widehat Q_-.$ 
Furthermore, $\mu = \lambda - \sum_{i = 0}^r k_i \alpha_i = - \alpha_{-1} - \sum_{i = 0}^r
k_i \alpha_i \in \widehat Q_-$, where the coefficients $k_i$ are as in 
(1.20).  
The weight multiplicity of $\mu$ in  $L(\lambda)$ is the same as the root
multiplicity mult$(\mu)$ of $\mu$ in ${\Cal L}$, which can be computed using
the following Freudenthal-type recursive formula due to Peterson. 
\vskip 2mm

\proclaim{ Proposition 3.1} \ {\rm ([P], cf. [K2, Exercise 11.12])} 
 \hskip 5mm For $\beta \in \widehat
Q_{-}$, define  

$$c_{\beta}=\sum_{n \ge 1} \dsize\frac {1}{n} \text {mult} 
\left(\dsize\frac {\beta}{n}\right).$$

\noindent Then 
$$(\beta|\beta+2\rho) c_{\beta}
=\sum \Sb \beta', \beta'' \in \widehat Q_{-} \\ \beta=\beta'+\beta'' \endSb
(\beta'|\beta'') c_{\beta'} c_{\beta''},  \tag 3.2$$ 

\noindent where $\rho \in {\frak h}^\ast$ is such that $\rho(h_i) = 1$ for
$i = -1, 0, \cdots, r$.   
\endproclaim

\vskip 2mm
\hskip 5mm We write $\mu=\mu_0-d_{\lambda}(\mu) \delta$, where 
$d_{\lambda}(\mu)$ and $\mu_0=-\alpha_{-1}-\sum_{i=0}^r m_i \alpha_i$ 
are given by (1.19) and (1.21). Since the coefficient 
of $\alpha_{-1}$ in $\mu$ is $-1$, any decomposition
is of the form $\mu=\beta'+\beta''$, where
$$\beta'=-\alpha_{-1}-\sum_{i=0}^r s_i \alpha_i, \ \ \text {and} \ 
\  \beta''=-\sum_{i=0}^r t_i \alpha_i $$
with $s_i, t_i \in \text {\bf Z}_{\ge 0}$ or it has the form
with the roles of $\beta'$ and $\beta''$ switched.
Note that $c_{\mu}=\text {mult} (\mu)$ and $c_{\beta'}=\text {mult}
(\beta')$, where mult$(\cdot)$ is the multiplicity in ${\Cal L}$, which, 
for $\mu$ and $\beta'$, is the same as the multiplicity in $L(\lambda)$. 
Thus, in order to have a nontrivial contribution to $c_{\beta'}$ and
$c_{\beta''}$, $\beta'$ must be a weight of $L(\lambda)$ 
and $\beta''=-k\alpha$
for some $k\ge 1$, where $\alpha$ is a positive root of ${\frak g}$.

\vskip 2mm

\proclaim{ Lemma 3.3} \ Suppose $d_{\lambda}(\mu) > 0$. Then
 $(\mu|\mu+2\rho)$ is a polynomial
in
$r$ of degree 1.
\endproclaim
\vskip 2mm
\noindent
{\it Proof.} \ Since 
$\mu=-\alpha_{-1}-\sum_{i=0}^r m_i \alpha_i -d_{\lambda}(\mu) \delta$, 
we have 
$$
\aligned
& (\mu|\mu+2\rho)=(-\alpha_{-1}-\sum_{i=0}^r m_i \alpha_i 
-d_{\lambda}(\mu) \delta | -\alpha_{-1}-\sum_{i=0}^r m_i \alpha_i 
-d_{\lambda}(\mu) \delta + 2\rho) \\
& \ \ = -2 d_{\lambda}(\mu) (r+1) -2 \sum_{i=0}^r m_i a_i
-2 \sum_{i=0}^r m_i +\sum_{i,j=0}^r m_im_j a_{i,j}
-2l d_{\lambda}(\mu)\\
& \quad \quad \quad  +(\alpha_{-1}|\alpha_{-1})-2(\rho|\alpha_{-1}).
\endaligned
$$
By (1.21), the terms 
$\sum_{i=0}^r m_i a_i$, $\sum_{i=0}^r m_i$, 
$\sum_{i,j=0}^r m_i m_j a_{i,j}$ are all constants. Therefore
$(\mu|\mu+2\rho)$ is a polynomial in $r$ of degree 1.
\hskip 1cm $\square$

\vskip 2mm
\hskip 5mm We now state and prove our main result.

\vskip 2mm
\proclaim{ Theorem 3.4} \ Let ${\frak g}$ be the affine Kac-Moody algebra 
of type $A_r^{(1)}$, and let $\lambda\in P^+$
be a dominant integral weight of level $l>0$ for ${\frak g}$
with determining data
$\underline a=(a_0, a_1, \cdots, a_{s-1})$,
$\underline a'=(a_{r-t+1}, a_{r-t+2}, \cdots, a_{r})$,
and $m \in \text {\bf Z}$.
Assume that $r \ge l+s+t+2$, and let $\mu \in P^+$ be a dominant 
integral weight of level $l$ with determining data
$\underline b=(b_0, b_1, \cdots, b_{s'-1})$,
$\underline b'=(b_{r-t'+1}, b_{r-t'+2}, \cdots, b_{r})$,
and $n \in \text {\bf Z}$
such that $s'+t' \le r,\;  s+t' \le r$, and $s'+t \le r$.
Suppose that $\mu$ is related to $\lambda$ for infinitely many values of $r\ge l+s+t+2$.
If $\mu$ is a weight of $L(\lambda)$ for some $r_{0} \ge l+s+t+2$, then
it is a weight of $L(\lambda)$ for all $r \ge r_0$, and 
the multiplicity of $\mu$ in $L(\lambda)$ is given by a polynomial in $r$ of degree
$\le d_{\lambda}(\mu)$.
\endproclaim

\vskip 2mm
\noindent    
{\it Proof.} \ We will prove our assertion by induction on $d_{\lambda}(\mu)$
and on the partial ordering on the affine weight lattice. 
Write 
$\mu=-\alpha_{-1}-\sum_{i=0}^r m_i \alpha_i -d_{\lambda}(\mu) \delta$,
where $d_{\lambda}(\mu)$ and the $m_i$'s are given by (1.19) and (1.21).

\vskip 2mm
\hskip 5mm If $d_{\lambda}(\mu)<0$, then for $i=p-1, p, \cdots, r-q, r-q+1$,
the coefficient of $\alpha_{i}$ in $\mu$ is positive. Hence $\mu$ cannot 
be a weight of $L(\lambda)$, and therefore its multiplicity in $L(\lambda)$
is zero. 
\vskip 2mm
\hskip 5mm  Let $p=\text {max} (s, s')$, $q=\text {max} (t,t')$.
Suppose that $d_{\lambda}(\mu)=0$. Then 
$$\mu=\lambda-\sum_{i=0}^{p-2} m_i \alpha_i
-\sum_{i=r-q+2}^r m_i \alpha_i.$$

The multiplicity of $\mu$ in $L(\lambda)$ is the
number of linearly independent vectors of the form
$f_{i_1} f_{i_2} \cdots f_{i_k} \cdot v_{0}$, where
$v_0$ is the highest weight vector of $L(\lambda)$ and 
$f_j$ appears $m_j$ times in the expression for each $j \in
\{0, 1, \cdots, p-2, r-q+2, r-q+3, \cdots, r \}$. Clearly, this number is
independent of $r$ (it may be 0).
In particular, if this number is nonzero
for some $r_{0}\ge l+s+t+2$, 
then it is nonzero and constant for all $r\ge r_0$.
Therefore, if $\mu$ is a weight of $L(\lambda)$ for some $r_0\ge l+s+t+2$,
then it is a weight of $L(\lambda)$ for all $r\ge r_0$, and 
the multiplicity of $\mu$ is a constant. 
The same argument shows that
any $\tau \in P$ such that $d_{\lambda} (\tau)=0$ 
has a constant multiplicity (which may be 0) in $L(\lambda)$. 

\vskip 2mm
\hskip 5mm Suppose $d_{\lambda}(\mu) \ge 1$. Consider a 
dominant integral weight $\tau$ of level $l$ with determining data 
$\underline c=(c_0, c_1, \cdots, c_{s''-1})$,
$\underline c'=(c_{r-t''+1}, c_{r-t''+2}, \cdots, c_{r})$,
and $n' \in \text {\bf Z}$
such that $s''+t'' \le r, \; s+t'' \le r$, and $s''+t \le r$
which is related to $\lambda$ for infinitely many values of 
$r\ge l+s+t+2$, and write
$\tau=\tau_0-d_{\lambda}(\tau) \delta$,
where $d_{\lambda}(\tau_0)=0$.
Assume that if $d_{\lambda}(\tau) <d_{\lambda}(\mu)$ 
or if $d_{\lambda}(\tau)=d_{\lambda}(\mu)$ and $\tau_0 > \mu_0$,
our assertion holds for $\tau$. That is, we assume that
if $\tau$ is a weight of $L(\lambda)$ for some $r_0\ge l+s+t+2$,
then it is a weight of $L(\lambda)$ for all $r\ge r_0$, and
the multiplicity of $\tau$ in $L(\lambda)$ is a polynomial
in $r$ of degree $\le d_{\lambda}(\tau)$.

\vskip 2mm
\hskip 5mm Consider a decomposition $\mu=\beta'+\beta''$,
where $\beta' \in -\alpha_{-1}-Q_- = \lambda - Q_-$ and $\beta''$
is a  multiple of a negative root of ${\frak g}$. Thus we may assume that
$\beta''$ is one of the following: 
$$ \text {(i)} \  -k\delta, \ \ \text {(ii)} \  -k\gamma,  \ \ 
\text {(iii)} \  -k( k' \delta+\gamma), \ \ \text {(iv)} \ 
-k( k'\delta-\gamma), \tag 3.5$$ 
where $k, k' \ge 1$ and $\gamma$ is a positive root of ${\frak g}_0=A_r$. 

\vskip 2mm
\hskip 5mm Note that $\beta''=-k\delta$ is an imaginary root of ${\frak g}$
for all $k\ge 1$, and its multiplicity in ${\frak g}$ (and hence in
${\Cal L}$) is $r$ (see [K1, Cor. 7.4]). It follows that
$$
\aligned
c_{-k\delta}&=\sum_{m\ge 1} \dsize\frac {\text {mult} (-k\delta /m)}{m} 
=\sum_{m|k} \dsize\frac {\text {mult} ( -(k/m) \delta)}{m} \\
&= \sum_{m|k} \dsize\frac {r}{m}=\dsize\frac {\xi(k)}{k} r,
\endaligned
\tag 3.6   
$$  
where $\xi(k)$ denotes the sum of all factors of $k$.
On the other hand, the roots $\gamma$, $k' \delta \pm \gamma$ are real, and 
their multiplicities in ${\frak g}$ (and hence in ${\Cal L}$) are all 1.
Moreover, if $k\ge 2$, the multiplicities of $k\gamma$ and $k( k' \delta \pm
\gamma)$ in ${\frak g}$ are all 0. 
Therefore we have
$$c_{\beta''}=\dsize\frac {1}{k} \ \ \text {if} \ 
\beta''=-k\gamma \ \text {or} \ \beta''=-k(k' \delta
\pm \gamma). \tag 3.7$$

\vskip 2mm
\hskip 5mm We now treat the four cases separately.

\vskip 2mm  
\noindent 
{\bf Case 1.} \ Suppose first that $\beta''=-k\delta$ for $k\ge 1$. In this
case, 
$$\beta'=-\alpha_{-1}-\sum_{i=0}^r m_i \alpha_i 
-(d_{\lambda}(\mu)-k) \delta \in \widehat Q_{-}.$$
Since $m_i=0$  for $i=p-1, p, \cdots, r-q, r-q+1$, we must have $d_{\lambda}(\mu)-k \ge 0$ in
order for $\beta'$ to belong to $\widehat Q_{-}$, which implies that $k$ runs from 1 to 
$d_{\lambda}(\mu)$.  Now
$$
\aligned
 (\beta'|\beta'')&=(-\alpha_{-1}-\sum_{i=0}^r m_i \alpha_i 
-(d_{\lambda}(\mu)-k) \delta | -k\delta)\\
&= k(\alpha_{-1}|\delta)=-kl,
\endaligned
$$
where $l$ denotes the level. 
As we have seen before, if $\beta'$ is not a weight of $L(\lambda)$,
then $c_{\beta'}=0$ and there is no contribution to the
right-hand side of (3.2).   So we may assume that $\beta'$ is a weight of 
$L(\lambda)$ for some $r_0 \ge l+s+t+2$.   Observe that  $\beta'$ is
dominant since $\beta'(h_j) = \mu(h_j) \geq 0$ for all $j \in I$.   
Since $\beta'$ is related to $\lambda$ for 
infinitely many values of $r\ge l+s+t+2$, and since
$d_{\lambda}(\mu)-k < d_{\lambda}(\mu)$, it follows from the
induction hypothesis that $\beta'$ is a weight of $L(\lambda)$ 
for all $r\ge r_0$, and 
$\text {mult} (\beta')$ is a polynomial in $r$ of degree 
$\le d_{\lambda}(\mu)-k \le d_{\lambda}(\mu)-1.$
We have seen in (3.6) that $c_{\beta''}=r \xi(k)/k$, which is a polynomial
in $r$ of degree 1. Therefore, the total contribution of the various
decompositions of this kind to the right-hand side of (3.2) is a polynomial in
$r$ of degree $\le d_{\lambda}(\mu)$.    

\vskip 2mm
\noindent
{\bf Case 2.} \ Suppose $\beta''=-k\gamma$ for $k\ge 1$, where $\gamma$ is 
a positive root of ${\frak g}_0=A_r$. 
Thus $\gamma=\alpha_u+\alpha_{u+1}+\cdots + \alpha_{v}$ with $1\le u \le v \le r$. 
In this case, we have 
$$\beta'=-\alpha_{-1}-\sum_{i=0}^r m_i \alpha_i +k\gamma 
-d_{\lambda}(\mu) \delta \in \widehat Q_{-}.$$
Note that for all $i=u,u+1, \cdots, v$, the coefficient of $\alpha_i$
in $\beta'$ must be $\le 0$. That is, $-m_i+k-d_{\lambda}(\mu) \le 0$
for $i=u,u+1, \cdots, v$.
Let $M=\text {max} \{m_i | \ 0 \le i \le r \}$. 
Then $k\le M+d_{\lambda}(\mu)$, and hence $k$ ranges from 1 to $M+
d_{\lambda}(\mu)$. 
Note that $M$ is independent of the value of $r$. 
We also have
$$
\aligned
(\beta'|\beta'')&=(-\alpha_{-1}-\sum_{i=0}^r m_i \alpha_i +k\gamma 
-d_{\lambda}(\mu) \delta | -k\gamma)\\
&= k(\alpha_{-1}|\gamma)+k \sum_{i = 0}^r m_i (\alpha_i|\gamma)
-k^2 (\gamma|\gamma).
\endaligned
$$
Hence $(\beta'|\beta'')$ is a constant for each $k$,  
because $(\lambda|\alpha_i)=-a_i=0$ and $m_i=0$ for 
$i=p,p+1, \cdots, r-q$,

\vskip 2mm
\hskip 5mm 
Now by (1.5), 
$\gamma=\alpha_u+\alpha_{u+1}+\cdots+\alpha_v 
= -\Lambda_{u-1}+\Lambda_u+ \Lambda_v-\Lambda_{v+1}$. 
Observe that 
$${\frak r}_{u-1}(-\Lambda_{u-1}+\Lambda_u+\Lambda_v-\Lambda_{v+1})
=-\Lambda_{u-2}+\Lambda_{u-1}+\Lambda_v-\Lambda_{v+1}, \quad
\text{and}$$
 $${\frak r}_{v+1}(-\Lambda_{u-1}+\Lambda_u+\Lambda_v-\Lambda_{v+1})
=-\Lambda_{u-1}+\Lambda_{u}+\Lambda_{v+1}-\Lambda_{v+2},$$ 

\noindent where ${\frak r}_{i}$ denotes the simple
reflection corresponding to the root
$\alpha_i$. Therefore for each $r\ge l+s+t+2$, 
we apply the simple reflections 
${\frak r}_{u-1}$, ${\frak r}_{u-2}$, 
$\cdots$, ${\frak r}_1$, ${\frak r}_0$, ${\frak r}_r$, $\cdots$
and then ${\frak r}_{v+1}$, ${\frak r}_{v+2}$, $\cdots,$ ${\frak r}_{r}$,
${\frak r}_{0}$, ${\frak r}_{1}$, $\cdots$ in succession
to get a dominant integral weight. 
(It may take several rounds of applying the simple reflections 
in this order to produce a dominant integral weight.)
Let $w_{r}$ denote the corresponding Weyl group element of $A_{r}^{(1)}$.
We can verify that $w_{r} \beta'$ has the form
$$w_{r} \beta'=\sum_{i=0}^r b'_i (r) \Lambda_i-n'(r)\delta,$$
where $b'_i (r) =0$ for $i=p+1,\cdots,r-q-1.$
Moreover, it is tedious but straightforward to show that 
the sequences of integers $\underline {c} (r)
= (b_{0}'(r), b_{1}'(r), \cdots, b_{p}'(r);$
$b_{r-q}'(r), b_{r-q+1}'(r), \cdots, b_{r}'(r); n'(r))$ 
are the same for all $r\ge l+s+t+2$. 
That is, the determining data of $w_{r} \beta'$ is given by 
$\underline c=(b'_0, b'_1, \cdots, b'_{x})$,
$\underline c'=(b'_{r-y}, b'_{r-y+1}, \cdots, b'_{r})$,
and $n' \in \text {\bf Z}$,
where $x\le p$, $y \le q$ and 
$b_{i}'=b_{i}'(r)$, $n'=n'(r)$ \ (for all $r\ge l+s+t+2$).
Rather than writing $w_{r} \beta'$ in what
follows, we denote the dominant weight determined by this
data as $\tau$.  

Since $\tau$ is related to $\lambda$, 
$\tau$ also has level $l$, and 
if we let $d'_i=b'_i-a_i$ for $i=0,1,\cdots,r$,
then $a_i=b'_i=d'_i=0$ for $i=p+1,p+2,\cdots,r-q-1$. Hence we may write
$$\tau=-\alpha_{-1}-\sum_{i=0} ^r m'_i\alpha_i
-d_{\lambda}(\tau)\delta,$$
where $-\alpha_{-1}-\sum_{i=0} ^r m'_i\alpha_i$ has depth 0 with respect to $\lambda$.
In addition, since $\tau$ is the highest element among the
Weyl group conjugates of $\beta'$, the inequality   
$\mu< \beta' \le \tau$ must hold, and hence 
$$\tau-\mu =\sum_{i=0}^r (m_i-m_i') \alpha_i
+(d_{\lambda}(\mu)-d_{\lambda}(\tau))\delta \in \widehat Q_{+},
\tag 3.8$$
and
$$\tau-\beta'=\sum_{i=0}^r (m_i-m_i') \alpha_i -k\gamma 
+ (d_{\lambda}(\mu)-d_{\lambda}(\tau)) \delta \in \widehat Q_{+}.
\tag 3.9
$$
By the same argument as in Lemma 1.23, we can show that $m'_i \ge 0$ for all
$i=p-1,p,\cdots,r-q,r-q+1$, and hence
$d_{\lambda}(\tau)\le d_{\lambda}(\mu)$.  

\vskip 2mm
\hskip 5mm If $d_{\lambda}(\tau)=d_{\lambda}(\mu)$, then  
$m_i' \le 0$ by (1.26), and hence $m_i'=0$ for $i=p-1, p, \cdots, r-q, r-q+1$.
Since $m_i=0$ for $i=p-1, p, \cdots, r-q, r-q+1$, (3.9) can be written as
$$\tau-\beta'=\sum_{i=0}^{p-2} (m_i-m_i') \alpha_i
+\sum_{i=r-q+2}^r (m_i-m_i')\alpha_i -k\gamma \in \widehat Q_{+}.
\tag 3.10
$$
Thus, in order for $\tau-\beta'$ to be an element of $\widehat Q_{+}$,
$\gamma$ must be a linear combination of the simple roots 
$\alpha_1, \alpha_2, \cdots, \alpha_{p-2}$, or of 
$\alpha_{r-q+2}, \cdots, \alpha_{r-1},  \alpha_r$.
The number of such $\gamma$ is at most 
$\displaystyle{\frac {(p-2)(p-1)}{2}+\frac {(q-1)q}{2}}$, which is independent 
of $r$.
Since $\tau > \mu$ and $d_{\lambda}(\tau)=d_{\lambda}(\mu)$,
we have $\tau_0 > \mu_0$.  Hence by the induction hypothesis, 
if $w_{r_{0}}\beta'$ is a weight of $L(\lambda)$ for some $r_{0}\ge 
l+s+t+2$, then $\tau = w_{r}\beta'$ is a weight of $L(\lambda)$ for all 
$r\ge r_0$, and 
$\text {mult} (\beta')=\text {mult} (\tau)$ is a polynomial in $r$ 
of degree $\le d_{\lambda}(\mu)$.
Note that $c_{\beta''}=1/k$ for all $k= 1, \cdots, M+d_{\lambda}(\mu)$. 
Therefore, the contribution of 
these partitions to the right-hand side of (3.2) is a polynomial in 
$r$ of degree $\le d_{\lambda}(\mu)$. 

\vskip 2mm
\hskip 5mm Suppose that $d_{\lambda}(\tau) < d_{\lambda}(\mu)$.
By the induction hypothesis, if $w_{r_{0}} \beta'$ is a weight 
of $L(\lambda)$ for some $r_{0} \ge l+s+t+2$, then 
$\tau = w_{r}\beta'$ is a weight of $L(\lambda)$ for all
$r\ge r_0$, and  $\text {mult} (\beta')=\text {mult}(\tau)$ 
is a polynomial of degree
$\le d_{\lambda}(\tau) \le d_{\lambda}(\mu)-1$. 
Since there are $\displaystyle{\frac {r(r+1)}{2}}$ positive roots in ${\frak g}_0=A_r$,
a polynomial in $r$  of degree 2, and since $c_{\beta''}=1/k$ for all
$k=1,\cdots, M+d_{\lambda}(\mu)$, the 
contribution of these decompositions to the right-hand side of (3.2) is a 
polynomial in $r$ of degree $\le d_{\lambda}(\mu)+1.$ 

\vskip 2mm
\hskip 5mm Therefore, the total
contribution of the partitions in Case 2 to the right-hand side of (3.2)
is a polynomial in $r$ of degree $\le d_{\lambda}(\mu)+1$. 

\vskip 2mm
\noindent
{\bf Case 3.} \ Suppose that $\beta''=-k( k' \delta+\gamma)$, where $k,k'
\ge 1$ and $\gamma$ is a positive root of ${\frak g}_0=A_r$. 
In this case, 
$$
\aligned
\beta' &= -\alpha_{-1}-\sum_{i=0}^r m_i \alpha_i +k( k' \delta+\gamma)
-d_{\lambda}(\mu) \delta \\
&= -\alpha_{-1}-\sum_{i=0}^r m_i \alpha_i +k\gamma
-(d_{\lambda}(\mu) -k k' ) \delta \in \widehat Q_{-}.
\endaligned
$$
Observe that the coefficient of $\alpha_0$ in $\beta'$ is
$-m_0-d_{\lambda}(\mu) +k k'$, which must be $\le 0$. Thus
$k k' \le m_0 + d_{\lambda}(\mu)$, and hence $k, k'$ range from 1 to
$m_0+d_{\lambda}(\mu)$. 
We have 
$$
\aligned
(\beta'|\beta'')&=( -\alpha_{-1}-\sum_{i=0}^r m_i \alpha_i +k\gamma
-(d_{\lambda}(\mu) -k k' ) \delta | -k(k'\delta + \gamma)) \\
&= -k k' l + k(\alpha_{-1}|\gamma)+k\sum_{i=0}^r m_i (\alpha_i|\gamma)
-k^2 (\gamma|\gamma),
\endaligned
$$
which can be seen to be a constant as in Case 2. 

\vskip 2mm
\hskip 5mm Moreover, by the same argument as in Case 2, for 
each $r\ge l+s+t+2$, we can verify that 
$\beta'$ is Weyl group conjugate to a dominant 
integral weight $\tau = w_{r} \beta'$ that has the form
$\tau=\sum_{i=0}^r b'_i\Lambda_i-n'\delta,$
where $b'_i=0$ for $i=p+1,\cdots,r-q-1$.
So if we let 
$d_i'=b_i'-a_i$ $(i=0,1,\cdots,r)$, we may write
$$\tau=-\alpha_{-1}-\sum_{i=0}^r m_i' \alpha_i 
-d_{\lambda}(\tau) \delta, $$
where $-\alpha_{-1}-\sum_{i=0}^r m_i' \alpha_i$ has depth 0 
with respect to $\lambda$.
In addition, we have
$$\tau-\mu =\sum_{i=0}^r (m_i-m_i') \alpha_i
+(d_{\lambda}(\mu)-d_{\lambda}(\tau))\delta \in \widehat Q_{+},
\tag 3.11$$
and
$$
\aligned
\tau-\beta'&=\sum_{i=0}^r (m_i-m_i') \alpha_i -k\gamma \\
& \quad \quad  + (d_{\lambda}(\mu)-d_{\lambda}(\tau)-k k') \delta
\in \widehat Q_{+}.
\endaligned
\tag 3.12
$$
By the same argument as in Lemma 1.23,
we can show that $m_i'\ge 0$ for all $i=p-1,p, \cdots, r-q, r-q+1$, and 
hence $d_{\lambda}(\mu)-d_{\lambda}(\tau)-k k' \ge 0$. Therefore,
$$d_{\lambda}(\tau) \le d_{\lambda}(\mu) -k k' \le 
d_{\lambda}(\mu)-1.$$
By the induction hypothesis, 
if $w_{r_0} \beta'$ is a weight of $L(\Lambda)$ for some $r_{0}
\ge l+s+t+2$, then 
$\tau = w_{r}
\beta'$ is a weight of $L(\lambda)$ for all $r\ge r_0$, and 
$\text {mult} (\beta')=\text {mult} (\tau)$
is a polynomial in $r$ of degree $\le d_{\lambda}(\tau)
\le d_{\lambda}(\mu)-1$. Since there are $\displaystyle{\frac {r(r+1)}{2}}$
positive roots in ${\frak g}_0$ and since $c_{\beta''}=1/k$ for all $k=1,
\cdots, m_0+d_{\lambda}(\mu)$,
the total contribution of the partitions in this case to the right side of
(3.2) is a polynomial in $r$ of degree $\le d_{\lambda}(\mu)+1$.

\vskip 2mm
\noindent
{\bf Case 4.} \ Suppose $\beta''=-k(k'\delta-\gamma)$, where $k,k' \ge 1$ 
and $\gamma$ is a positive root of ${\frak g}_0=A_r$. In this case,
$$\beta'=-\alpha_{-1}-\sum_{i=0}^r m_i \alpha_i -k\gamma 
-(d_{\lambda}(\mu)-k k') \delta \in \widehat Q_{-}.$$
As in Case 3,  by looking at the coefficient of $\alpha_0$
in the above expression, we can show that 
$k$ and $k'$ range from 1 to $m_0+d_{\lambda}(\mu)$.
We have
$$
\aligned
(\beta'|\beta'')&= (-\alpha_{-1}-\sum_{i=0}^r m_i \alpha_i -k\gamma 
-(d_{\lambda}(\mu)-k k') \delta | -k(k' \delta -\gamma))\\
&= -k k' l-k(\alpha_{-1}| \gamma) -\sum_{i=0}^r m_i (\alpha_i|\gamma)
-k^2 (\gamma|\gamma),
\endaligned
$$
which can be seen to be a constant as in Case 2. 

\vskip 2mm
\hskip 5mm  By the identical argument as in Case 2
we can verify for each $r\ge l+s+t+2$ that $\beta'$ is Weyl group conjugate to a dominant
integral weight $\tau = w_{r}\beta'$ that has the form
$\tau =\sum_{i=0}^r b'_i\Lambda_i-n'\delta,$
where $b'_i=0$ for $i=p+1,\cdots,r-q-1$.
So if we suppose as before 
$d_i'=b_i'-a_i$ $(i=0,1,\cdots,r)$, then  
$$\tau =-\alpha_{-1}-\sum_{i=0}^r m_i' \alpha_i
-d_{\lambda}(\tau) \delta,$$
where $-\alpha_{-1}-\sum_{i=0}^r m_i' \alpha_i $ has depth 0
with respect to $\lambda$. 
Moreover,  
$$\tau-\mu =\sum_{i=0}^r (m_i-m_i') \alpha_i
+(d_{\lambda}(\mu)-d_{\lambda}(\tau))\delta \in \widehat Q_{+},
\tag 3.13$$
and
$$
\aligned
\tau-\beta' & =\sum_{i=0}^r (m_i-m_i') \alpha_i + k\gamma \\
& \hskip 5mm + (d_{\lambda}(\mu)-d_{\lambda}(\tau)-k k') \delta
\in \widehat Q_{+}.
\endaligned
\tag 3.14
$$
Let us write $\gamma=\alpha_u +\alpha_{u+1} + \cdots + \alpha_{v}$
with $1\le u \le v \le r$. Then (3.14) becomes
$$
\aligned
\tau-\beta'&=\sum_{i=0}^r (m_i-m_i') \alpha_i-k(\alpha_0+\alpha_1+
\cdots + \alpha_{u-1}) \\
&\quad -k(\alpha_{v+1} +\cdots + \alpha_r)+ 
(d_{\lambda}(\mu)-d_{\lambda}(\tau)+k-k k') \delta
\in \widehat Q_{+}.
\endaligned
\tag 3.15
$$
The argument in Lemma 1.23 proves that $m_i'\ge 0$
for all $i=p-1,p, \cdots, r-q, r-q+1$ and that 
$d_{\lambda}(\mu)-d_{\lambda}(\tau)+k-k k' \ge 0$,  
which yields
$$d_{\lambda}(\tau) \le d_{\lambda}(\mu) +k-k k' \le d_{\lambda}(\mu).$$

\vskip 2mm
\hskip 5mm If $d_{\lambda}(\tau)=d_{\lambda}(\mu)$, then we must have
$k'=1$ and (3.15) can be written as 
$$
\aligned
\tau -\beta' & =\sum_{i=0}^r (m_i-m_i') \alpha_i \\
& \hskip 3mm -k(\alpha_0+\alpha_1+\cdots + \alpha_{u-1})
-k(\alpha_{v+1} +\cdots + \alpha_r) \in \widehat Q_{+}.
\endaligned
\tag 3.16
$$
Recall that $m_i=0$ and $m_i'\ge 0$ for $i=p-1,p, \cdots, r-q, r-q+1$.
Thus in order for $\tau-\beta'$ to belong to 
$\widehat Q_{+}$, it must be
that $m_i'=0$ for all $i=p-1,p, \cdots, r-q, r-q+1$, and 
$u\le p-2$, $v\ge r-q+1$. Hence the number of such $\gamma$ is at most 
$(p-2)(q-1)$, a constant. 
As $\tau>\mu$ and 
$d_{\lambda}(\tau)= d_{\lambda}(\mu)$ hold, 
we have $\tau_0 > \mu_0$.  Hence, 
it follows from the induction hypothesis that 
if $w_{r_{0}}\beta'$ is a weight of $L(\lambda)$ for some $r_{0} \ge l+s+t+2$, 
$\tau = w_{r}\beta'$ is a weight of $L(\lambda)$ for all $r\ge r_0$, and 
$\text {mult}(\beta')=\text {mult} 
(\tau)$ is a polynomial in $r$ of degree 
$\le d_{\lambda}(\tau)=
 d_{\lambda}(\mu)$. 
Note that $c_{\beta''}=1/k$ for all $i=1, \cdots, m_0+d_{\lambda}(\mu)$.
Therefore, the contribution of
these decompositions to the right-hand side of (3.2) is a polynomial in
$r$ of degree $\le d_{\lambda}(\mu)$.

\vskip 2mm
\hskip 5mm Suppose that $d_{\lambda}(\tau) < d_{\lambda}(\mu)$.
Then by the induction hypothesis,
if $w_{r_0}\beta'$ is a weight of $L(\lambda)$ for some $r_{0} \ge l+s+t+2$, 
then $\tau$ is a weight of $L(\lambda)$ for all $r\ge r_0$, and 
$\text {mult} (\beta')=\text {mult}
(\tau)$ is a polynomial in $r$ of degree $\le d_{\lambda}
(\tau) \le d_{\lambda}(\mu)-1$. 
Since there are $\displaystyle{\frac {r(r+1)}{2}}$ positive roots in 
${\frak g}_0=A_r$ and since $c_{\beta''}=1/k$ for all $k=1,\cdots, 
m_0+d_{\lambda}(\mu)$, the contribution of these 
decompositions to the right side of (3.2) is a polynomial in $r$ of degree 
$\le d_{\lambda}(\mu)+1$. 

\vskip 2mm
\hskip 5mm Therefore, what the   
partitions in Case 4 contribute to the right-hand side of (3.2) is a polynomial 
in $r$ of degree $\le d_{\lambda}(\mu)+1$.

\vskip 2mm
\hskip 5mm Consequently, the sum of all the contributions from
Case 1 to Case 4, which is the right side of (3.2),  is a polynomial 
in $r$ of degree $\le d_{\lambda}(\mu)+1$. 
By Lemma 3.3 and (3.2), 
we have $\text {mult}(\mu)=f/g$, where $f$ is a polynomial
in $r$ of degree $\le d_{\lambda}(\mu)+1$ and $g$ is a polynomial
in $r$ of degree 1.
Since $\text {mult}(\mu)$ takes positive integral values for 
infinitely many values of 
$r\ge l+s+t+2$,  it must be a polynomial in $r$ (see [PS], p. 130), and 
$$\text {deg} (\text {mult} (\mu)) \le 
(d_{\lambda}(\mu)+1)-1=d_{\lambda}(\mu).$$
This completes the proof of the theorem. \hskip 1cm $\square$

\vskip 8mm
\example{Example 3.17}
The following tables illustrate the polynomial behavior of the
multiplicity of the weight $\mu-k\delta$ 
in the irreducible highest weight module $L(\lambda)$ over the affine
Kac-Moody algebra $A_r^{(1)}$.  
The numerical data in these tables was taken
from [KMPS]. From now on, let  $\underline a=(a_0, a_1, \cdots, a_{s-1})$,
$\underline a'=(a_{r-t+1}, a_{r-t+2}, \cdots, a_{r})$,
and $m \in \text {\bf Z}$ be the determining data for $\lambda$, and let
 $\underline b=(b_0, b_1, \cdots, b_{s'-1})$,
$\underline b'=(b_{r-t'+1}, b_{r-t'+2}, \cdots, b_{r})$,
and $n \in \text {\bf Z}$ be the determining data for $\mu.$      
\endexample
\vfill \eject 
\noindent 1. \ \ \ $\lambda=\mu=\Lambda_0+\Lambda_r;$

\quad $\underline a=(1)$, $\underline a'=(1)$, $m=0,$ \; \;
$\underline b=(1)$, $\underline b'=(1)$, $n=0$;

\quad $d_{\lambda}(\mu-k\delta)=(k-m)-(d_1+2d_2+\cdots+(p-1)d_{p-1})=k.$  

\bigskip

{\offinterlineskip
 \tabskip=0pt
 \halign{ \strut \vrule#& \, # \,  & \vrule#& \, # \, & \vrule#& \, # \,  & \vrule#& \, # \,  & 
\vrule#& \, # \, & \vrule#& \, # \,  & \vrule#& \, # \, & \vrule#& \, # \,  &  \vrule#& \, # \, &  
\vrule#& \, # \, & \vrule#  \cr
\noalign{\hrule} 
&&&&&&&&&&&&&&&&&&&& \cr
\noalign{\vskip -10pt}
& k $\setminus$ r  && 1 && 2 && 3 && 4 && 5 && 6 && 7 && 8 && polynomial & \cr  
&&&&&&&&&&&&&&&&&&&& \cr
\noalign{\vskip -10pt}
\noalign{\hrule} 
&&&&&&&&&&&&&&&&&&&& \cr
\noalign{\vskip -10pt}
& 0 && 1 && 1 && 1 && 1 && 1 && 1 && 1 && 1 && 1 & \cr  
&&&&&&&&&&&&&&&&&&&& \cr
\noalign{\vskip -10pt}
&&&&&&&&&&&&&&&&&&&& \cr
\noalign{\vskip -10pt}
& 1 && 2 && 4 && 6 && 8 && 10 && 12 && 14 && 16 && $2r$ & \cr  
&&&&&&&&&&&&&&&&&&&& \cr
\noalign{\vskip -10pt}
&&&&&&&&&&&&&&&&&&&& \cr
\noalign{\vskip -10pt}
& 2 && 4 && 13 && 27 && 46 && 70 && 99 && 133 && 172 && $\frac{1}{2}(5r^2+3r)$ \hfil & \cr  
&&&&&&&&&&&&&&&&&&&& \cr
\noalign{\vskip -10pt}
&&&&&&&&&&&&&&&&&&&& \cr
\noalign{\vskip -10pt}
& 3 && 8 && 36 && 98 && 208 && 380 && 628 && 966 && 1408 && $\frac{1}{3}(7r^3+9r^2$
\hfil & \cr
&&&&&&&&&&&&&&&&&&&& \cr
\noalign{\vskip -10pt} 
& && && && && && && && && && \quad \ \ $+8r)$ \hfil & \cr 
&&&&&&&&&&&&&&&&&&&& \cr
\noalign{\vskip -10pt}
& 4 && 14 && 89 && 310 && 804 && 1740 && 3329 && 5824 && 9520 && 
$\frac{1}{12}(21r^4+44r^3$ \hfil& \cr
&&&&&&&&&&&&&&&&&&&& \cr
\noalign{\vskip -10pt}  
& && && && && && && && && && \quad \  $+87r^2 +16r)$ \hfil & \cr
&&&&&&&&&&&&&&&&&&&& \cr
\noalign{\vskip -10pt}
&&&&&&&&&&&&&&&&&&&& \cr
\noalign{\vskip -10pt}
& 5 && 24 && 204 && 888 && 2768 && 7012 && 15396 && 30436 && 55520 && 
$\frac{1}{30}(33r^5+100r^4$ \hfil & \cr
&&&&&&&&&&&&&&&&&&&& \cr
\noalign{\vskip -10pt} 
& && && && && && && && && && \quad \ \  $+315r^3$ \hfil & \cr
&&&&&&&&&&&&&&&&&&&& \cr
\noalign{\vskip -10pt} 
& && && && && && && && && && \quad \quad $+200r^2$ \hfil & \cr
&&&&&&&&&&&&&&&&&&&& \cr
\noalign{\vskip -10pt} 
& && && && && && && && && && \quad \quad \ \  $+72r)$ \hfil & \cr
\noalign{\vskip -10pt} 
&&&&&&&&&&&&&&&&&&&& \cr
\noalign{\hrule}
}}

\bigskip
\bigskip
\bigskip
\bigskip 
\noindent 2. \ \ \  $\lambda=3\Lambda_1, \quad  \mu=\Lambda_0+\Lambda_1+\Lambda_2;$

\quad $\underline a=(0,3)$, $\underline a'=\emptyset$, $m=0,$ \; \;
            $\underline b=(1,1,1)$, $\underline b'=\emptyset$, $n=0;$

\quad $d_{\lambda}(\mu-k\delta)=(k-m)-(d_1+2d_2+\cdots+(p-1)d_{p-1})=k.$ 

\bigskip

{\offinterlineskip
 \tabskip=0pt
 \halign{ \strut \vrule#& \, # \,  & \vrule#& \, # \,  & 
\vrule#& \, # \, & \vrule#& \, # \,  & \vrule#& \, # \, & \vrule#& \, # \,  &  \vrule#& \, # \, &  \vrule#& \, # \, & \vrule#  \cr
\noalign{\hrule} 
&&&&&&&&&&&&&&&&\cr
\noalign{\vskip -10pt}
& k $\setminus$ r  &&  2 && 3 && 4 && 5 && 6 && 7 &&  polynomial & \cr  
&&&&&&&&&&&&&&&& \cr
\noalign{\vskip -10pt}
\noalign{\hrule} 
&&&&&&&&&&&&&&&& \cr
\noalign{\vskip -10pt}
& 0 && 1 && 1 && 1 && 1 && 1 && 1 && 1 & \cr  
&&&&&&&&&&&&&&&& \cr
\noalign{\vskip -10pt}
&&&&&&&&&&&&&&&& \cr
\noalign{\vskip -10pt}
& 1 &&  4 && 6 && 8 && 10 && 12 && 14 && $2r$ & \cr  
&&&&&&&&&&&&&&&& \cr
\noalign{\vskip -10pt}
&&&&&&&&&&&&&&&& \cr
\noalign{\vskip -10pt}
& 2 && 15 && 31 && 53 && 81 && 115 && 155 &&  $3r^2+r+1$ \hfil & \cr  
&&&&&&&&&&&&&&&& \cr
\noalign{\vskip -10pt}
&&&&&&&&&&&&&&&& \cr
\noalign{\vskip -10pt}
& 3 && 44 && 126 && 278 && 523 && 884 && 1384 && $\frac{1}{6}(23r^3+3r^2+40r-12)$ \hfil & \cr
&&&&&&&&&&&&&&&& \cr
\noalign{\vskip -10pt}
&&&&&&&&&&&&&&&& \cr
\noalign{\vskip -10pt}
& 4 && 121 && 456 && 1267 && 2901 && 5808 && 10541 && $\frac{1}{24}(103r^4-54r^3+533r^2$\hfil &
\cr   
&&&&&&&&&&&&&&&& \cr
\noalign{\vskip -10pt}
& && && && && && && && \quad \ \  $-294r+144)$ \hfil & \cr
&&&&&&&&&&&&&&&& \cr
\noalign{\vskip -10pt}
&&&&&&&&&&&&&&&& \cr
\noalign{\vskip -10pt}
& 5 && 300 && 1477 && 5120 && 14166 && 33444 && 70188 && $\frac{1}{120}(513r^5-800r^4$
\hfil  & \cr
&&&&&&&&&&&&&&&& \cr
\noalign{\vskip -10pt}  
& && && && && && && && \quad \ \  $+5815r^3-6580r^2$ \hfil& \cr
&&&&&&&&&&&&&&&& \cr
\noalign{\vskip -10pt}  
& && && && && && && && \quad \quad   $+7412r-22)$ \hfil & \cr 
\noalign{\vskip -10pt}
&&&&&&&&&&&&&&&& \cr
\noalign{\hrule}
}}

\vfill \eject
\noindent 3. \ \ \  $\lambda=2\Lambda_0, \quad  \mu=\Lambda_1+\Lambda_r$;
   
\quad $\underline a=(2)$, $\underline a'=\emptyset$, 
$m=0,$ \; \; $\underline b=(1)$, $\underline b'=(1)$, $n=0$;

\quad $d_{\lambda}(\mu-k\delta)=(k-m)-(d_1+2d_2+\cdots+(p-1)d_{p-1})=k-1$.  

\bigskip

{\offinterlineskip
 \tabskip=0pt
 \halign{ \strut  \vrule#& \, # \, & \vrule#& \, # \,  & \vrule#& \, # \,  & 
\vrule#& \, # \, & \vrule#& \, # \,  & \vrule#& \, # \, & \vrule#& \, # \,  
&  \vrule#& \, # \, &  \vrule#& \, # \, & \vrule#  \cr
\noalign{\hrule} 
&&&&&&&&&&&&&&&&&& \cr
\noalign{\vskip -10pt}
& k $\setminus$ r  && 2 && 3 && 4 && 5 && 6 && 7 && 8 && polynomial & \cr  
&&&&&&&&&&&&&&&&&& \cr
\noalign{\vskip -10pt}
\noalign{\hrule} 
&&&&&&&&&&&&&&&&&& \cr
\noalign{\vskip -10pt}
& 0 && 0 && 0 && 0 && 0 && 0 && 0 && 0 && 0 & \cr  
&&&&&&&&&&&&&&&&&& \cr
\noalign{\vskip -10pt}
&&&&&&&&&&&&&&&&&& \cr
\noalign{\vskip -10pt}
& 1 && 1 && 1 && 1 && 1 && 1 && 1 && 1 && 1&  \cr  
&&&&&&&&&&&&&&&&&& \cr
\noalign{\vskip -10pt}
&&&&&&&&&&&&&&&&&& \cr
\noalign{\vskip -10pt}
& 2 && 4 && 6 && 8 && 10 && 12 && 14 && 16 &&  $2r$ & \cr  
&&&&&&&&&&&&&&&&&& \cr
\noalign{\vskip -10pt}
&&&&&&&&&&&&&&&&&& \cr
\noalign{\vskip -10pt}
& 3 && 12 && 25 && 43 && 66 && 94 && 127 && 165  && $\frac{1}{2}(5r^2+r+2)$ & \cr
&&&&&&&&&&&&&&&&&& \cr
\noalign{\vskip -10pt}
&&&&&&&&&&&&&&&&&& \cr
\noalign{\vskip -10pt}
& 4 && 32 && 87 && 186 && 343 && 572 && 887 && 1302  && $\frac{1}{3}(7r^3+3r^2+17r-6)$ \hfil &
\cr   
&&&&&&&&&&&&&&&&&& \cr
\noalign{\vskip -10pt}
& 5 && 77 && 266 && 693 && 1513 && 2923 && 5162 && 8511 && 
$\frac{1}{12}(21r^4+16r^3+129r^2$ & \cr
&&&&&&&&&&&&&&&&&& \cr
\noalign{\vskip -10pt}  
& && && && && && && && && \quad \ \ $-46r+36)$ \hfil \cr
\noalign{\vskip -10pt}
&&&&&&&&&&&&&&&&&& \cr
\noalign{\hrule}
}}

\bigskip
\bigskip
\bigskip
\bigskip 
\noindent 4. \ \ \  $\lambda=\Lambda_0+\Lambda_1, \quad  \mu=\Lambda_2+\Lambda_r;$

\quad $\underline a=(1,1)$, $\underline a'=\emptyset$, 
$m=0,$ \; \; $\underline b=(0,0,1)$, $\underline b'=(1)$, $n=0;$

\quad $d_{\lambda}(\mu-k\delta)=(k-m)-(d_1+2d_2+\cdots+(p-1)d_{p-1})=k-1.$ 
\bigskip

{\offinterlineskip
 \tabskip=0pt
 \halign{ \strut \vrule#& \, # \,  & \vrule#& \, # \,  & 
\vrule#& \, # \, & \vrule#& \, # \,  & \vrule#& \, # \, & \vrule#& \, # \,  &  \vrule#& \, # \, &  \vrule#& \, # \, & \vrule#  \cr
\noalign{\hrule} 
&&&&&&&&&&&&&&&&\cr
\noalign{\vskip -10pt}
& k $\setminus$ r  &&  3 && 4 && 5 && 6 && 7 && 8 &&  polynomial & \cr  
&&&&&&&&&&&&&&&& \cr
\noalign{\vskip -10pt}
\noalign{\hrule} 
&&&&&&&&&&&&&&&& \cr
\noalign{\vskip -10pt}
& 0 && 0 && 0 && 0 && 0 && 0 && 0 && 0 &  \cr  
&&&&&&&&&&&&&&&& \cr
\noalign{\vskip -10pt}
&&&&&&&&&&&&&&&& \cr
\noalign{\vskip -10pt}
& 1 &&  2 && 2 && 2 && 2 && 2 && 2 && 2  & \cr  
&&&&&&&&&&&&&&&& \cr
\noalign{\vskip -10pt}
&&&&&&&&&&&&&&&& \cr
\noalign{\vskip -10pt}
& 2 && 12 && 17 && 22 && 27 && 32 && 37 &&  $5r-3$ \hfil & \cr  
&&&&&&&&&&&&&&&& \cr
\noalign{\vskip -10pt}
&&&&&&&&&&&&&&&& \cr
\noalign{\vskip -10pt}
& 3 && 50 && 92 && 148 && 218 && 302 && 400 && $7r^2-7r+8$ \hfil & \cr
&&&&&&&&&&&&&&&& \cr
\noalign{\vskip -10pt}
&&&&&&&&&&&&&&&& \cr
\noalign{\vskip -10pt}
& 4 && 172 && 396 && 770 && 1336 && 2136 && 3212 && $7r^3-9r^2+28r-20$ \hfil & \cr  
&&&&&&&&&&&&&&&& \cr
\noalign{\vskip -10pt}
&&&&&&&&&&&&&&&& \cr
\noalign{\vskip -10pt}
& 5 && 522 && 1466 && 3382 && 6816 && 12446 && 21082 && $\frac{1}{2}(11r^4-16r^3+97r^2$ \hfil  
& \cr &&&&&&&&&&&&&&&& \cr
\noalign{\vskip -10pt}  
& && && && && && && && \quad \ \ $-124r+84)$ \hfil \cr
\noalign{\vskip -10pt}
&&&&&&&&&&&&&&&& \cr
\noalign{\hrule}
}}

\vfill \eject

\noindent 5. \ \ \ $\lambda=2\Lambda_0+\Lambda_1, \,\, \mu=2\Lambda_1+\Lambda_r;$

\quad $\underline a=(2,1)$, $\underline a'=\emptyset$, 
$m=0,$ \;\; $\underline b=(0,2)$, $\underline b'=(1)$, $n=0;$

\quad $d_{\lambda}(\mu-k\delta)=(k-m)-(d_1+2d_2+\cdots+(p-1)d_{p-1})=k-1.$ 

\bigskip

{\offinterlineskip
 \tabskip=0pt
 \halign{ \strut \vrule#& \, # \,  & \vrule#& \, # \,  & 
\vrule#& \, # \, & \vrule#& \, # \,  & \vrule#& \, # \, & \vrule#& \, # \,  &  \vrule#& \, # \, &  \vrule#& \, # \, & \vrule#  \cr
\noalign{\hrule} 
&&&&&&&&&&&&&&&&\cr
\noalign{\vskip -10pt}
& k $\setminus$ r  && 2 &&  3 && 4 && 5 && 6 && 7 &&   polynomial & \cr  
&&&&&&&&&&&&&&&& \cr
\noalign{\vskip -10pt}
\noalign{\hrule} 
&&&&&&&&&&&&&&&& \cr
\noalign{\vskip -10pt}
& 0 && 0 && 0 && 0 && 0 && 0 && 0 && 0 & \cr  
&&&&&&&&&&&&&&&& \cr
\noalign{\vskip -10pt}
&&&&&&&&&&&&&&&& \cr
\noalign{\vskip -10pt}
& 1 &&  1 && 1 && 1 && 1 && 1 && 1 && 1 & \cr  
&&&&&&&&&&&&&&&& \cr
\noalign{\vskip -10pt}
&&&&&&&&&&&&&&&& \cr
\noalign{\vskip -10pt}
& 2 && 6 && 9 && 12 && 15 && 18 && 21 &&  $3r$ & \cr  
&&&&&&&&&&&&&&&& \cr
\noalign{\vskip -10pt}
&&&&&&&&&&&&&&&& \cr
\noalign{\vskip -10pt}
& 3 && 22 && 49 && 87 && 136 && 196 && 267 && $\frac{1}{2}(11r^2-r+2)$ \hfil & \cr
&&&&&&&&&&&&&&&& \cr
\noalign{\vskip -10pt}
&&&&&&&&&&&&&&&& \cr
\noalign{\vskip -10pt}
& 4 && 70 && 214 && 492 && 951 && 1638 && 2600 && $\frac{1}{6}(47r^3-21r^2+76r-24)$ \hfil & \cr 
&&&&&&&&&&&&&&&& \cr
\noalign{\vskip -10pt}
&&&&&&&&&&&&&&&& \cr
\noalign{\vskip -10pt}
& 5 && 193 && 795 && 2328 && 5515 && 11304 && 20868 && $\frac{1}{8}(75r^4-86r^3+373r^2$ \hfil &
\cr 
&&&&&&&&&&&&&&&& \cr
\noalign{\vskip -10pt}  
& && && && && && && && \quad \ \  $-290r+120)$ \hfil & \cr
\noalign{\vskip -10pt}
&&&&&&&&&&&&&&&& \cr
\noalign{\hrule}
}}

\bigskip
\bigskip
\vskip 8mm
\baselineskip=12pt

\centerline {\bf References \rm}
\vskip 4mm 
 
\item{[BKM1]} G. Benkart, S.-J. Kang,  K. C. Misra, 
Graded Lie algebras of Kac-Moody type, Adv. in Math. {\bf 97} (1993), 154-190.
\vskip 2.2mm

\item{[BKM2]} G. Benkart, S.-J. Kang, K. C. Misra, 
Weight multiplicity polynomials for affine Kac-Moody algebras of 
type $A_r^{(1)}$, Compositio Math. {\bf 104} (1996), 153-187. 
\vskip 2.2mm

\item{[BK]} G. Benkart, S. N. Kass, 
Weight multiplicities for affine Kac-Moody algebras, in {\it Modern Trends in Lie
Theory}, Queen's Papers in Pure and Applied Math. {\bf 94}, V. Futorny and D.
Pollack eds. (1994), 1-12. \vskip 2.2mm

\item{[FF]} A. J. Feingold, I. B. Frenkel, A hyperbolic Lie algebra
and the theory of Siegel modular forms of genus 2,  Math. Ann. {\bf 263}
(1983), 87-144. 
  
\vskip 2.2mm
\item{[GK]} O. Gabber, V. G. Kac, On defining relations of certain
infinite-dimensional  Lie algebras, Bull. Amer. Math. Soc. {\bf 5} (1981),
185-189.  \vskip 2.2mm

\item{[K1]} V. G. Kac, Simple irreducible graded Lie algebras of finite
growth, Math. USSR-Izv. {\bf 2} (1968), 1271-1311.
\vskip 2.2mm

\item{[K2]} V. G. Kac, {\it Infinite Dimensional Lie
Algebras}, 3rd ed.,  Cambridge Univ. Press, Cambridge,  1990.
\vskip 2.2mm

\item{[Ka1]} S.-J. Kang, {\it Gradations and Structure of Kac-Moody Lie 
Algebras}, Yale University Ph.D. dissertation  1990.
\vskip 2.2mm 

\item{[Ka2]} S.-J. Kang, {\it Root multiplicities of Kac-Moody  
algebras}, Duke Math. J. {\bf 74} (1994), 635-666.
\vskip 2.2mm 

\item{[KMPS]} S. N. Kass, R. V. Moody, J. Patera, R. Slansky,
{\it Affine Lie Algebras, Weight Multiplicities, and Branching Rules, 
Vols. I and II}, Los
Alamos Series in Basic and Applied Sciences, Univ. of Calif. Press,  1990.
\vskip 2.2mm

\item{[P]} D. H. Peterson, Freudenthal-type formulas for root and
weight multiplicities, preprint (unpublished).
\vskip 2.2mm

\item{[PS]} G. P\'{o}lya and G. Szeg\"{o}, 
{\it Problems and Theorems in Analysis II}, Springer-Verlag,  1976.

\enddocument